\documentclass{article}
\usepackage{amssymb,amsmath}
\usepackage{graphics}

\def\qed{\hfill {\hbox{${\vcenter{\vbox{               
   \hrule height 0.4pt\hbox{\vrule width 0.4pt height 6pt
   \kern5pt\vrule width 0.4pt}\hrule height 0.4pt}}}$}}}

\def\tr{\triangleright}

\newtheorem{theorem}{Theorem}
\newtheorem{definition}{Definition}

\newtheorem{proposition}[theorem]{Proposition}

\newtheorem{example}{Example}
\newtheorem{conjecture}{Conjecture}

\newenvironment{proof}[1][Proof]{\smallskip\noindent{\bf #1.}\quad}%
{\qed\par\medskip}

\date{}

\title{\Large \textbf{A polynomial invariant of finite quandles}}

\author{Sam Nelson \\ 
\small Department of Mathematics, Whittier College \\
\small 13406 Philadelphia,  P.O. Box 634 \\
\small Whittier, CA 90608-0634 \\
\texttt{knots@esotericka.org}}

\begin{document}
\maketitle

\begin{abstract}
We define a two-variable polynomial invariant of finite quandles. 
In many cases this invariant completely determines the algebraic structure 
of the quandle up to isomorphism. We use this polynomial to define
a family of link invariants which generalize the quandle counting invariant.
\end{abstract}



\section{\large \textbf{Introduction}}

A \textit{quandle} is an algebraic structure whose axioms are transcriptions
of the three Reidemeister moves, making quandles a natural source of
invariants of knots and links. In \cite{J}, Joyce associated a finitely
generated quandle to each tame knot or link  $L\subset S^3$ and proved that 
this knot quandle determines the knot type up to orientation-reversing 
homeomorphism of pairs $(L,S^3)$.

Since then, much work has been done in the study of knot invariants defined
in terms of the knot quandle, such as  \cite{D2}, \cite{C4}, etc. In 
particular, for any finite quandle $T$ there are invariants of knots and
links defined in terms of the set $\mathrm{Hom}(Q(K),T)$ of quandle 
homomorphisms from the knot quandle $Q(K)$ into a finite target quandle $T$. 
Since isomorphic target quandles $T$ and $T'$ define the same invariants, we 
would like to have a convenient way of distinguishing finite quandles. 
Moreover, understanding the structure of these finite quandles can lead to 
a better understanding of their link invariants.

In this paper we describe a two-variable polynomial invariant of finite 
quandles. This polynomial encodes a set with multiplicities arising from 
counting trivial actions of elements on other elements; the polynomial 
is a convenient way to encode this information with the added benefit of 
obtaining integer-valued invariants via specialization of $s$ and $t$.
While not a complete invariant, this polynomial does seem to be 
quite good at distinguishing finite quandles. Indeed, in some cases the
invariant completely determines the structure of the quandle.

One motivation for this construction comes from consideration of the nature of
identity elements and kernels of homomorphisms. Specifically, groups have a 
unique identity element which acts trivially on the rest of the group, and 
the preimage of this element under a homomorphism measures how much 
information about the domain group is lost in the homomorphism. Unlike groups,
quandles do not have a unique identity element; rather, every element in a 
quandle acts trivially on itself and possibly some other elements. We might 
say that in a group, the ``trivial action'' is concentrated in a single 
element, while in a quandle it is distributed throughout the set. Our 
polynomial may then be understood as a way of quantifying how this ``trivial 
action'' is distributed.

The paper is organized as follows. In section 2, we recall the definition
of a quandle and some examples of quandles. We then define the $qp$ polynomial
and list some examples. In section 3, we define a generalization of $qp$
for subquandles of a larger quandle and a $qp$-based Laurent polynomial
associated to a homomorphism between finite quandles. We the use $qp$ to 
define a family of link invariants which are jazzed-up versions of the 
familiar quandle counting invariant.
In section 4, we list some questions for future research and end with a 
table listing the $qp$ values for all quandles with order up to 5. 

\section{\large \textbf{Quandles and a polynomial invariant}}

We begin with a definition from Joyce, \cite{J}.

\begin{definition} \textup{ A \textit{quandle} is a set
$Q$ with a binary operation $\tr:Q\times Q\to Q$ satisfying
\newcounter{c}
\begin{list}{(\roman{c})}{\usecounter{c}} 
\item{for every $a\in Q$, $a\tr a=a$,}
\item{for every $a,b\in Q$, there is a unique $c\in Q$ such that $a=c\tr b$, 
and}
\item{for every $a,b,c\in Q$, $(a\tr b)\tr c=(a\tr c)\tr(b\tr c).$}
\end{list}
If (ii) and (iii) are satisfied, then $Q$ is a \textit{rack} \cite{FR};
if (iii) is satisfied, then $Q$ is a \textit{shelf} \cite{CC}.}
\end{definition}

The axioms can be understood as transcriptions of the three 
Reidemeister moves where each arc in the knot diagram corresponds to a quandle
element and the quandle operation corresponds to an arc going through a 
crossing. Some well-known examples of quandles include:

\begin{example}
\textup{A group $G$ is a quandle under various operations, such as
\newcounter{d}
\begin{list}{\arabic{d}.}{\usecounter{d}}
\item{$x\tr y = y^{-1} xy$, called a \textit{conjugation quandle},}
\item{$x\tr y = y^{-n}xy^n$, an $n$-\textit{fold conjugation quandle},}
\item{$x\tr y = s(xy^{-1})y$ where $s\in \mathrm{Aut(G)}$, called a 
\textit{homogeneous quandle}.}
\end{list}}
\end{example}

\begin{example}
\textup{A commutative ring $R$ with a choice of invertible element $t\in R$
is a quandle under}
\[x\tr y = tx+(1-t) y.\] \textup{Indeed, it's easy to see that this is 
\textit{only} possible quandle structure on $R$ with the quandle operation
defined as a linear polynomial: set $x\tr y = Ax+By+C.$ Then axiom (i)
requires 
\[x=Ax+Bx+C = (A+B)x + C\] for all $x\in R$, so comparing coefficients
we have $C=0$ and $B=1-A$. More generally, we may replace the ring $R$ with 
a module over $\mathbb{Z}[t^{\pm 1}]$; quandles of this type are called 
\textit{Alexander quandles}. See \cite{N} for more.}
\end{example}

\begin{example} \textup{Let $R$ be a commutative ring, $M$ an $R$-module and 
$\langle,\rangle:M\times M\to R$ an anti-symmetric bilinear form. Then $M$ is
a quandle under}
\[x\tr y = x+\langle x,y\rangle y.\]
\textup{If $R$ is a field and $\langle,\rangle$ is non-degenerate, then we
call $M$ a \textit{symplectic quandle}. See \cite{Y} and \cite{NN} for more.}
\end{example}

Let $Q=\{x_1, x_2, \dots, x_n\}$ be a finite quandle. We can describe the
algebraic structure of $Q$ with an $n\times n$ matrix $M_Q$ satisfying
$M_Q[i,j]=k$ where $x_i\tr x_j=x_k$. In other words, the matrix of $Q$ is
just the quandle operation table without the ``$x$''s. This quandle matrix
notation is convenient for symbolic computation purposes \cite{HN}, \cite{HMN}.

\begin{definition}\label{qp}
\textup{Let $Q$ be a finite quandle. For any element $x\in Q$, let}
\[ c(x) = \left|\{ y\in Q \ : \ y\tr x = y\}\right|
\quad \mathrm{and \ let} 
\quad r(x) = \left|\{ y\in Q \ : \ x\tr y = x \}\right| 
.\]
\textup{Then we define the \textit{quandle polynomial of $Q$}, $qp_Q(s,t)$, 
to be $\displaystyle{qp_Q(s,t)=\sum_{x\in Q} s^{r(x)}t^{c(x)}.}$ We will 
refer to $r(x)$ and $c(x)$ as the \textit{row count} and \textit{column count}
of $x$ respectively.}
\end{definition}

That is, $r(x)$ is the number of elements of $Q$ which act trivially on $x$, 
while $c(x)$ is the number of elements of $Q$ on which $x$ acts trivially.
These numbers are easily computed from the matrix of $Q$, simply by going 
through rows and columns and counting occurrences of the row number. Note 
that if $Q$ has infinite cardinality, then $c(x)$ and $r(x)$ may not be finite
and $qp_Q(s,t)$ may not be defined.

\begin{proposition}\label{iso}
If $Q$ and $Q'$ are isomorphic finite quandles, then $qp_Q(s,t)=qp_{Q'}(s,t)$.
\end{proposition}

\begin{proof}
Suppose $f:Q\to Q'$ is an isomorphism of quandles and fix $x\in Q$. 
Then for every $y\in C(x)=\{ y\in Q \ : \ y\tr x = y\}$, we have
$f(y)\tr f(x) = f(y\tr x) =f(y)$, so $f(y)\in C(f(x))$ and we have
$|C(x)|\le |C(f(x)|$. Repeating the same argument with $f^{-1}$ yields the
opposite inequality, and we have $c(x)=c(f(x))$. A similar argument
shows that $r(x)=r(f(x))$. Then  we have
\[qp_Q(s,t)=\sum_{x\in Q} s^{r(x)}t^{c(x)}
=\sum_{f(x)\in Q'}s^{r(x)}t^{c(x)}
=\sum_{f(x)\in Q'}s^{r(f(x))}t^{c(f(x))}=qp_{Q'}(s,t).
\]
\end{proof}

\begin{definition}
\textup{The specialization $qp_Q(1,t)$ is the \textit{column polynomial} of
$Q$. The specialization $qp_Q(s,1)$ is the \textit{row polynomial} of $Q$.}
\end{definition}

\begin{example}
\textup{The \textit{trivial quandle of order $n$}, $T_n$, satisfies 
$x_i\tr x_j=x_i$ for all $i,j\in \{1,2,\dots, n\}.$ It has quandle matrix}
\[M_{T_n}=\left[
\begin{array}{cccc}
1 & 1 & \dots & 1 \\
2 & 2 & \dots & 2 \\
\vdots & \vdots & & \vdots \\
n & n & \dots & n
\end{array}
\right]
\quad \mathrm{and} \quad qp_{T_n}(s,t) = ns^nt^n.
\]
\end{example}

\begin{proposition}
Evaluating $qp_Q(1,1)$ yields $|Q|$, the cardinality of $Q$.
\end{proposition}

\begin{proof}
\[qp(1,1)=\sum_{x\in Q} 1^{r(x)}1^{c(x)} = \sum_{x\in Q} 1 =|Q|.\]
\end{proof}

Thus, we may regard $qp_Q(s,t)$ as a kind of enhanced version of cardinality 
for finite quandles. It's easy to see that $qp_Q(s,t)$ completely determines 
the quandle structure for some quandles, such as $T_n$ above. In fact, $qp$ 
distinguishes all quandles of order 3 and 4 and all non-Latin quandles of order 5,
as shown in table 1.

\begin{example}
\textup{There is only one possible quandle $Q$ with $qp_Q(s,t)=3st$. To see 
this, consider the matrix of a quandle $Q$ with $qp_Q(s,t)=3st$. This is a 
$3\times 3$ quandle
matrix such that every row and column of $M_Q$ has precisely one entry equal 
to its row number. We know (see \cite{HN}) that every quandle matrix must have
columns which are permutations of $\{1, 2, \dots, n\}$ and must have entries 
$1,2,\dots, n$ along the diagonal. Then we must have 
$M_Q=\left[\begin{array}{ccc}
1 & 0 & 0 \\
0 & 2 & 0 \\
0 & 0 & 3
\end{array}
\right]$
(here we use ``0'' as a blank to be filled in). Now, the $(2,1)$ entry cannot 
be $2$, or the element $x_1\in Q$ would have $c(x_1)>1$, and similarly the 
$(3,1)$ entry cannot be 3. Indeed, just the $c(x)$ numbers here are enough to 
determine that we must have  $M_Q=\left[\begin{array}{ccc}
1 & 3 & 2 \\
3 & 2 & 1 \\
2 & 1 & 3
\end{array}
\right]$.}
\end{example} 

There are examples of non-isomorphic quandles with the same $qp$ values,
as the next example shows.

\begin{example}
\textup{The Alexander quandles $\mathbb{Z}_5[t^{\pm 1}]/(t-2)$, 
$\mathbb{Z}_5[t^{\pm 1}]/(t-3)$ and $\mathbb{Z}_5[t^{\pm 1}]/(t-4)$ all have
$qp_Q(s,t)=5st$; the fact that they are non-isomorphic follows from corollary
2.2 in \cite{N}.}
\end{example}

The $qp$ polynomial distinguishes between some quandles which have 
the same orbit decomposition but different structure maps (see \cite{NW}).

\begin{example}
\textup{The quandles with quandle matrices}
\[
M_Q=\left[\begin{array}{cccccc}
1 & 3 & 2 & 1 & 1 & 1 \\
3 & 2 & 1 & 2 & 2 & 2 \\
2 & 1 & 3 & 3 & 3 & 3 \\
4 & 4 & 4 & 4 & 6 & 5 \\
5 & 5 & 5 & 6 & 5 & 4 \\
6 & 6 & 6 & 5 & 4 & 6  
\end{array}\right] \quad \mathrm{and} \quad
M_{Q'}=\left[\begin{array}{cccccc}
1 & 3 & 2 & 1 & 3 & 2 \\
3 & 2 & 1 & 3 & 2 & 1 \\
2 & 1 & 3 & 2 & 1 & 3 \\
4 & 6 & 5 & 4 & 6 & 5 \\
6 & 5 & 4 & 6 & 5 & 4 \\
5 & 4 & 6 & 5 & 4 & 6  
\end{array}\right]
\]
\textup{have $qp_Q(s,t)=6s^4t^4$ and $qp_{Q'}(s,t)=6s^2t^2$ respectively. Both 
have the same orbit decomposition, namely two orbits isomorphic to the 
Alexander quandle $\mathbb{Z}_3[t^{\pm 1}]/(t-2)$, but the two have different 
structure maps.}
\end{example}

\begin{definition}\textup{
A quandle is \textit{connected} if it has only one orbit. A quandle is
\textit{Latin} or \textit{strongly connected} if the quandle operation
is left-invertible, in addition to the right-invertibility required by 
axiom (ii).\footnote{Latin quandles are a type of distributive quasigroup.}
In particular, the rows of a Latin quandle are also permutations of 
$\{1, 2, \dots, n\}$}
\end{definition}

\begin{proposition}
A Latin quandle $Q$ satisfies $qp_Q(s,t)=|Q|st.$
\end{proposition}

\begin{proof}\label{lat}
Let $Q$ be a Latin quandle and fix an element $x_i$. Since $Q$ is a quandle, 
the $(i,i)$ entry of $M_Q$ is $i$, so $c(x_i)\ge 1$ and $r(x_i)\ge 1$. Now,
suppose $r(x_i)>1$. Then there is a column $j\ne i$ with $M_Q[i,j]=i$. But
then row $i$ is not a permutation of $\{1,2,\dots, n\}$, contradicting the
fact that $Q$ is Latin. Hence, we must have $r(x_i)=1$. 

Similarly, if $c(x_i)>1$ then there is some $j\ne i$ such that $M_Q[j,i]=j$,
and then we have $r(x_j)>1$, contradicting our conclusion above. 
Hence $c(x_i)=1$. Since $i$ was arbitrary, we have
\[qp_Q(s,t) = \sum_{x\in Q} st = |Q|st.\]


\end{proof}

Direct computations show that the converse of proposition \ref{lat} has no 
counterexamples with cardinality less than 7. Thus we have:

\begin{conjecture}
Every quandle $Q$ with $qp_Q(s,t)=|Q|st$ is Latin.
\end{conjecture}

We may also define the $qp$ invariant for non-quandle racks; such racks may 
have $qp_R(s,t)=0$, as in the next example.

\begin{example}\textup{
Let $R=\mathbb{Z}_n$ with $i\tr j = i+k$ for a fixed $k\in R$. Then $R$ is
a rack since}
\[ (x\tr y)\tr z = (x+k)\tr z = x+2k\] \textup{while}
\[ (x\tr z) \tr (y\tr z) =(x+k)\tr (y+k)=x+2k.\]\textup{
The matrix of $R$ has every entry in row $i$ equal to $i+k$, so if
$k\ne 0 \in \mathbb{Z}_n,$ then $c(x)=r(x)=0$ for all $x\in R$, 
and $qp_R(s,t)=0$.}
\end{example}

Indeed, if $qp_R(1,1)\ne|R|$ for a rack $R$, then $R$ is not a quandle.

\section{\large \textbf{Subquandles, homomorphisms and a family of
link invariants}}

In this section we define a generalization of $qp_Q(s,t)$ which gives 
information about how a subquandle is embedded in a larger quandle and we use 
this to define a family of link invariants.

\begin{definition}
\textup{Let $S\subset Q$ be a subquandle of $Q$. The \textit{subquandle
polynomial} $qp_{S\subset Q}(s,t)$ of $S\subset Q$ is}
\[qp_{S\subset Q}(s,t) =\sum_{x\in S}s^{r(x)}t^{c(x)}\]
\textup{where $r(x)$ and $c(x)$ are the same as in definition \ref{qp}.}
\end{definition}

\begin{example}
\textup{Let $Q$ be the quandle with matrix 
$\displaystyle{M_Q=\left[\begin{array}{cccc}
1 & 1 & 2 & 2 \\
2 & 2 & 1 & 1 \\
3 & 3 & 3 & 3 \\
4 & 4 & 4 & 4 
\end{array}\right]}$ and let $S=\{1,2\}$ and $S'=\{3,4\}$. Then both
$S$ and $S'$ are isomorphic to $T_2$, the trivial quandle of order 2, but they
are embedded in $S$ in a different way; this is reflected in their subquandle
polynomial values $qp_{S\subset Q}(s,t)=2s^2t^4$ and 
$qp_{S'\subset Q}(s,t)=2s^4t^2$.}
\end{example}

\begin{proposition}
If a finite quandle $Q=Q_1\cup Q_2 \cup \dots \cup Q_n$ is a disjoint 
union of several subquandles (e.g., the orbit subquandles) then we have
\[qp_Q(s,t)=\sum_{i=1}^{n} qp_{Q_i\subset Q}(s,t).\]
\end{proposition}

\begin{proof}
\end{proof}

\begin{definition}
\textup{Let $f:Q\to Q'$ be a homomorphism of finite quandles.
Let $K_{qp}(f)\in \mathbb{Z}[s^{\pm 1}, t^{\pm 1}]$ be given by}
\[K_{qp}(f)=\sum_{x\in Q} s^{r(f(x))-r(x)}t^{c(f(x))-c(x)}.\]
\end{definition}

\begin{proposition}\label{ker}
If $f:Q\to Q'$ is injective, then every exponent in $K_{qp}(f)$ is 
nonnegative; if $f$ is surjective, then every exponent in $K_{qp}(f)$ is
nonpositive. If $f$ is an isomorphism, then $K_{qp}(f)=|Q|$.
\end{proposition}

\begin{proof}
This is similar to the proof of proposition \ref{iso}.
\end{proof}

Thus, $K_{qp}(f)$ has a philosophical similarity to the kernel of a 
homomorphism. Note, however, that the converse of proposition \ref{ker} is 
not true, as the next example shows.

\begin{example}\textup{
Let $f:T_2\to T_3$ be the constant map $f(1)=f(2)=1\in T_3$ where $T_n$
is the trivial quandle of order $n$. Then we have $K_{qp}(f)=2s^{3-2}t^{3-2}
=2st$, though $f$ is not injective.}
\end{example}


Recall that for any knot or link $K$, there is an associated \textit{knot 
quandle} $Q(K)$, and that for a given finite quandle $T$ the set of
quandle homomorphisms 
\[ \mathrm{Hom}(Q(K),T)=\{ f:Q(K)\to T \ : \ f(x\tr y) =f(x)\tr f(y)\}\]
is a source of computable knot invariants. Specifically, we can take the
cardinality of the set, which gives us an integer-valued invariant. 
Alternatively, we can count the homomorphisms weighted by a 
cocycle in one of the various quandle cohomology theories (described in
\cite{C4}, \cite{C3} etc.); these cocycles provide a way of squeezing extra 
information out of the set of homomorphisms. We can use the subquandle
polynomials of the image of each homomorphism in a similar way.

\begin{definition}
\textup{Let $K$ be a link and $T$ a finite quandle. Then for every 
$f\in \mathrm{Hom}(Q(K),T)$, the image of $f$ is a subquandle of $T$.
Define the \textit{subquandle polynomial invariant} $\Phi_{qp}(K)$ to
be the set with multiplicities}
\[\Phi_{qp}(K) =\{
qp_{\mathrm{Im}(f)\subset T}(s,t) \ : \ f\in \mathrm{Hom}(Q(K),T)\}.\]
\end{definition}

Normally, we encode sets with multiplicities whose elements are integers as 
polynomials where the multiplicities appear as coefficients and the elements 
of the set are powers of a variable; in this case, however, the elements of the
set are already polynomials, and polynomials with polynomial powers seem
a little awkward. However, we can derive convenient polynomial-valued
specializations of $\Phi_{qp}(K)$ by choosing values of $s$ and $t$ in 
$\mathbb{Z}$ and evaluating.

\begin{definition}
\textup{Let $K$ be a link, $T$ a finite quandle, 
$s_0,t_0\in \mathbb{Z}$. Define the \textit{specialized subquandle polynomial 
invariant} $\Phi_{qp}(K,s_0,t_0)\in \mathbb{Z}[z^{\pm 1}]$ to
be the polynomial}
\[\Phi_{qp}(K,s_0,t_0) = \sum_{f\in \mathrm{Hom}(Q(K),T)}
z^{qp_{\mathrm{Im}(f)\subset T}(s_0,t_0)}.\]
\end{definition}

\begin{example}\textup{
If we specialize $s_0=t_0=0$, then we have }
\[\Phi_{qp}(K,0,0) = \sum_{f\in \mathrm{Hom}(Q(K),T)}
z^{0} =|\mathrm{Hom}(Q)K),T)|,\]
\textup{so the specializations of $\Phi_{qp}(K)$ are generalizations of the
quandle counting invariant.}
\end{example}

In general, specializations of $\Phi_{qp}(K)$ contain more information than the
unadorned counting invariant, as the next example shows.

\begin{example}\textup{
The links $L_1$ and $L_2$ have quandle counting invariant 
\[|\mathrm{Hom}(Q(L_1),T)|=|\mathrm{Hom}(Q(L_2),T)|=13\] with the
quandle $T$ with quandle matrix listed below. However, the specialized 
quandle polynomials
\[\Phi_{qp}(L_1,1,1)=5z+2z^2+6z^3\quad \mathrm{and} \quad 
\Phi_{qp}(L_2,1,1)=5z+2z^2+6z^5\]
distinguish the links.}
\[\includegraphics{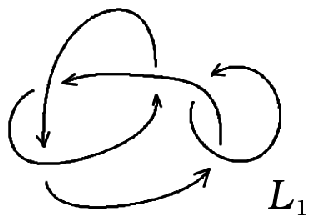} \quad \includegraphics{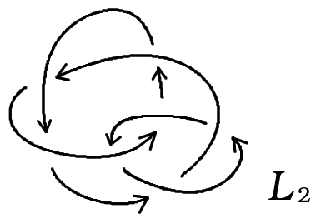}
\quad \raisebox{0.5in}{$M_T=\left[\begin{array}{ccccc}
1 & 3 & 2 & 3 & 2 \\
3 & 2 & 1 & 1 & 3 \\
2 & 1 & 3 & 2 & 1 \\
5 & 5 & 5 & 4 & 4 \\
4 & 4 & 4 & 5 & 5 
\end{array}\right].$}\]
\end{example}

\section{\large \textbf{Questions for further research}}

In this section we list a few questions whose answers may prove interesting.

Under what conditions are two finite quandles $qp$-equivalent? What 
information in addition to $qp_Q(s,t)$ must be specified in order to fully 
determine the quandle type of $Q$ up to isomorphism? For which polynomials
$p\in\mathbb{Z}[s,t]$ is the set of quandles $Q$ with $qp_Q(s,t)=p$ nonempty? 
It's easy to see that such a $p$ must have all positive coefficients whose sum
is $|Q|$ and 
each term must have degree $d$ satisfying  $2\le d\le 2|Q|$ (or $0\le d\le2|Q|$
in the rack case), but what other conditions are necessary or sufficient? 
Can we deduce from a quandle's polynomial whether the quandle is isomorphic
to a conjugation quandle, an Alexander quandle, a symplectic quandle, or a 
direct product, semidirect product, or an abelian extension of these?

For which specialized quandle polynomial invariants is there a quandle 
2-cocycle (see \cite{C4}) such that $\Phi_qp$ is equal to the quandle cocycle 
invariant? Do we gain any information by combining the two, perhaps by setting
\[\Psi_{\chi,qp}(K)=\sum_{f\in\mathrm{Hom}(Q(K),T)} 
z^{qp_{\mathrm{Im}(f)\subset T}(s_0,t_0)} t^{Co(f)}\]
where $Co(f)$ is the sum of the cocycle values at the crossings in $K$ 
(see \cite{C4})?

\begin{table}
\begin{center}
\begin{tabular}{|cc|cc|} \hline
& & & \\
$M_Q$ & $qp_Q(s,t)$ & $M_Q$ & $qp_Q(s,t)$ \\
& & & \\ \hline
& & & \\
$\left[\begin{array}{ccc} 
1 & 1 & 1 \\ 
2 & 2 & 2 \\ 
3 & 3 & 3 
\end{array}\right]$ 
& $3s^3t^3$ &
$\left[\begin{array}{ccc} 
1 & 1 & 1 \\ 
3 & 2 & 2 \\ 
2 & 3 & 3  
\end{array}\right]$ &
$2s^2t^3+s^3t$ \\
& & & \\
$\left[\begin{array}{ccc} 
1 & 3 & 2 \\ 
3 & 2 & 1 \\ 
2 & 1 & 3  
\end{array}\right]$ & $3st$ &
$\left[\begin{array}{cccc} 
1 & 1 & 1 & 1 \\ 
2 & 2 & 2 & 2 \\ 
3 & 3 & 3 & 3 \\ 
4 & 4 & 4 & 4 \\ 
\end{array}\right]$ & $4s^4t^4$ \\
& & & \\
$\left[\begin{array}{cccc} 
1 & 1 & 1 & 1 \\ 
2 & 2 & 2 & 3 \\ 
3 & 3 & 3 & 2 \\ 
4 & 4 & 4 & 4 \\ 
\end{array}\right]$ & $2s^3t^4+s^4t^2+s^4t^4$ &
$\left[\begin{array}{cccc} 
1 & 1 & 1 & 2 \\ 
2 & 2 & 2 & 3 \\ 
3 & 3 & 3 & 1 \\ 
4 & 4 & 4 & 4 \\ 
\end{array}\right]$ & $3s^3t^4+s^4t$ \\
& & & \\ 
$\left[\begin{array}{cccc} 
1 & 1 & 2 & 2 \\ 
2 & 2 & 1 & 1 \\ 
3 & 3 & 3 & 3 \\ 
4 & 4 & 4 & 4 \\ 
\end{array}\right]$ & $2s^2t^4+2s^4t^2$ &
$\left[\begin{array}{cccc} 
1 & 1 & 1 & 1 \\ 
2 & 2 & 4 & 3 \\ 
3 & 4 & 3 & 2 \\ 
4 & 3 & 2 & 4 \\ 
\end{array}\right]$ & $3s^2t^2+s^4t^4$ \\
& & & \\ 
$\left[\begin{array}{cccc} 
1 & 1 & 2 & 2 \\ 
2 & 2 & 1 & 1 \\ 
4 & 4 & 3 & 3 \\ 
3 & 3 & 4 & 4 \\ 
\end{array}\right]$ & $4s^2t^2$ &
$\left[\begin{array}{cccc} 
1 & 4 & 2 & 3 \\ 
3 & 2 & 4 & 1 \\ 
4 & 1 & 3 & 2 \\ 
2 & 3 & 1 & 4 \\ 
\end{array}\right]$ & $4st$ \\
& & & \\
$\left[\begin{array}{ccccc}
1 & 1 & 1 & 1 & 1 \\
2 & 2 & 2 & 2 & 2 \\
3 & 3 & 3 & 3 & 3 \\
4 & 4 & 4 & 4 & 4 \\
5 & 5 & 5 & 5 & 5
\end{array}\right]$ & $5s^5t^5$ &
$\left[\begin{array}{ccccc}
1 & 1 & 1 & 1 & 1 \\
2 & 2 & 2 & 2 & 3 \\
3 & 3 & 3 & 3 & 4 \\
4 & 4 & 4 & 4 & 2 \\
5 & 5 & 5 & 5 & 5
\end{array}\right]$ & $s^5t^2+3s^4t^5+s^5t^5$ \\
& & & \\
$\left[\begin{array}{ccccc}
1 & 1 & 1 & 1 & 2 \\
2 & 2 & 2 & 2 & 3 \\
3 & 3 & 3 & 3 & 4 \\
4 & 4 & 4 & 4 & 1 \\
5 & 5 & 5 & 5 & 5
\end{array}\right]$ & $4s^4t^5+s^5t$ &
$\left[\begin{array}{ccccc}
1 & 1 & 1 & 2 & 2 \\
2 & 2 & 2 & 3 & 3 \\
3 & 3 & 3 & 1 & 1 \\
4 & 4 & 4 & 4 & 4 \\
5 & 5 & 5 & 5 & 5
\end{array}\right]$ & $3s^3t^5+2s^5t^2$ \\
& & & \\
$\left[\begin{array}{ccccc}
1 & 1 & 1 & 1 & 1 \\
2 & 2 & 2 & 2 & 2 \\
3 & 3 & 3 & 5 & 4 \\
4 & 4 & 5 & 4 & 3 \\
5 & 5 & 4 & 3 & 5
\end{array}\right]$ & $3s^3t^3+2s^5t^5$ &
$\left[\begin{array}{ccccc}
1 & 1 & 2 & 2 & 2 \\
2 & 2 & 1 & 1 & 1 \\
3 & 3 & 3 & 3 & 3 \\
4 & 4 & 4 & 4 & 4 \\
5 & 5 & 5 & 5 & 5
\end{array}\right]$ & $2s^2t^5+3s^5t^3$ \\
& & & \\ \hline
\end{tabular}
\end{center}
\caption{$qp_Q(s,t)$ for quandles of order $\le 5$ part I.}
\end{table}

\begin{table}
\begin{center}
\begin{tabular}{|cc|cc|} \hline
& & & \\
$M_Q$ & $qp_Q(s,t)$ & $M_Q$ & $qp_Q(s,t)$ \\
& & & \\ \hline
& & & \\
$\left[\begin{array}{ccccc}
1 & 1 & 1 & 1 & 1 \\
2 & 2 & 2 & 2 & 2 \\
3 & 3 & 3 & 3 & 4 \\
4 & 4 & 4 & 4 & 3 \\
5 & 5 & 5 & 5 & 5
\end{array}\right]$ & $2s^4t^5+2s^5t^5+s^5t^3$ &
$\left[\begin{array}{ccccc}
1 & 1 & 1 & 1 & 2 \\
2 & 2 & 2 & 2 & 1 \\
3 & 3 & 3 & 3 & 4 \\
4 & 4 & 4 & 4 & 3 \\
5 & 5 & 5 & 5 & 5
\end{array}\right]$ & $4s^4t^5+s^5t$ \\
& & & \\
$\left[\begin{array}{ccccc}
1 & 1 & 1 & 1 & 1 \\
2 & 2 & 2 & 3 & 3 \\
3 & 3 & 3 & 2 & 2 \\
4 & 4 & 4 & 4 & 4 \\
5 & 5 & 5 & 5 & 5
\end{array}\right]$ & $2s^3t^5+s^5t^5+2s^5t^3$ &
$\left[\begin{array}{ccccc}
1 & 1 & 1 & 2 & 3 \\
2 & 2 & 2 & 3 & 1 \\
3 & 3 & 3 & 1 & 2 \\
4 & 4 & 4 & 4 & 4 \\
5 & 5 & 5 & 5 & 5
\end{array}\right]$ & $3s^3t^5+2s^5t^2$ \\
& & & \\
$\left[\begin{array}{ccccc}
1 & 1 & 1 & 2 & 2 \\
2 & 2 & 2 & 1 & 1 \\
3 & 3 & 3 & 3 & 3 \\
4 & 4 & 5 & 4 & 4 \\
5 & 5 & 4 & 5 & 5
\end{array}\right]$ & $2s^3t^5+2s^4t^3+s^5t^3$ &
$\left[\begin{array}{ccccc}
1 & 1 & 2 & 2 & 2 \\
2 & 2 & 1 & 1 & 1 \\
3 & 3 & 3 & 3 & 4 \\
4 & 4 & 4 & 4 & 3 \\
5 & 5 & 5 & 5 & 5
\end{array}\right]$ & $2s^2t^5+2s^4t^3+s^5t$ \\
& & & \\
$\left[\begin{array}{ccccc}
1 & 1 & 2 & 2 & 2 \\
2 & 2 & 1 & 1 & 1 \\
3 & 3 & 3 & 5 & 4 \\
4 & 4 & 5 & 4 & 3 \\
5 & 5 & 4 & 3 & 5
\end{array}\right]$ & $2s^2t^5+3s^3t$ &
$\left[\begin{array}{ccccc}
1 & 1 & 1 & 1 & 1 \\
2 & 2 & 5 & 3 & 4 \\
3 & 4 & 3 & 5 & 2 \\
4 & 5 & 2 & 4 & 3 \\
5 & 3 & 4 & 2 & 5
\end{array}\right]$ & $4s^2t^2+s^5t^5$ \\
& & & \\
$\left[\begin{array}{ccccc}
1 & 1 & 1 & 1 & 1 \\
2 & 2 & 2 & 3 & 3 \\
3 & 3 & 3 & 2 & 2 \\
4 & 5 & 5 & 4 & 4 \\
5 & 4 & 4 & 5 & 5
\end{array}\right]$ & $4s^3t^3+s^5t^5$ &
$\left[\begin{array}{ccccc}
1 & 1 & 2 & 2 & 2 \\
2 & 2 & 1 & 1 & 1 \\
3 & 3 & 3 & 3 & 3 \\
5 & 5 & 5 & 4 & 4 \\
4 & 4 & 4 & 5 & 5
\end{array}\right]$ & $4s^2t^3+s^5t$ \\
& & & \\
$\left[\begin{array}{ccccc}
1 & 1 & 1 & 1 & 1 \\
2 & 2 & 2 & 3 & 3 \\
3 & 3 & 3 & 2 & 2 \\
5 & 5 & 5 & 4 & 4 \\
4 & 4 & 4 & 5 & 5
\end{array}\right]$ & $2s^2t^3+2s^3t^3+s^5t^3$ &
$\left[\begin{array}{ccccc}
1 & 3 & 4 & 5 & 2 \\
3 & 2 & 5 & 1 & 4 \\
4 & 5 & 3 & 2 & 1 \\
5 & 1 & 2 & 4 & 3 \\
2 & 4 & 1 & 3 & 5
\end{array}\right]$ & $5st$ \\
& & & \\
$\left[\begin{array}{ccccc}
1 & 4 & 5 & 3 & 2 \\
3 & 2 & 4 & 5 & 1 \\
2 & 5 & 3 & 1 & 4 \\
5 & 1 & 2 & 4 & 3 \\
4 & 3 & 1 & 2 & 5
\end{array}\right]$ & $5st$ &
$\left[\begin{array}{ccccc}
1 & 4 & 5 & 2 & 3 \\
3 & 2 & 1 & 5 & 4 \\
4 & 5 & 3 & 1 & 2 \\
5 & 3 & 2 & 4 & 1 \\
2 & 1 & 4 & 3 & 5
\end{array}\right]$ & $5st$ \\
& & & \\ \hline
\end{tabular}
\end{center}
\caption{$qp_Q(s,t)$ for quandles of order $\le 5$  part II.}
\end{table}


\begin{thebibliography}{00}

\bibitem{CC}{J. S. Carter, A. Crans, M. Elhamdadi and M. Saito.
Cohomology of Categorical Self-Distributivity.
arXiv.org:math.GT/0607417}

\bibitem{C3}{J. S. Carter, M. Elhamdadi, M. Gra\~{n}a, and M. Saito.
Cocycle knot invariants from quandle modules and generalized quandle homology. 
\textit{Osaka J. Math.} \textbf{42} (2005) 499-541.}
  
\bibitem{C4}{J. S. Carter, D. Jelsovsky, S. Kamada, L. Langford and M. Saito.  
 Quandle cohomology and state-sum invariants of knotted curves and 
 surfaces.  \textit{Trans. Amer. Math. Soc.}  \textbf{355}  
(2003) 3947-3989.}

\bibitem{D2}{F. M. Dion\'{i}sio and P. Lopes. 
 Quandles at finite temperatures. II. 
 \textit{J. Knot Theory Ramifications} \textbf{12} (2003) 1041-1092.}

\bibitem{FR}{R. Fenn and C. Rourke.
 Racks and links in codimension two.
 \textit{J. Knot Theory Ramifications}  \textbf{1}  (1992), 343-406.}

\bibitem{HMN}{R. Henderson, T. Macedo and S. Nelson. 
Symbolic computation with finite quandles.
\textit{J. Symbolic Comput.} \textbf{41} (2006) 811-817.}

\bibitem{HN}{B. Ho and S. Nelson. Matrices and finite quandles.
\textit{Homology Homotopy Appl.} \textbf{7} (2005) 197-208.}

\bibitem{J}{D. Joyce.
 A classifying invariant of knots, the knot quandle.
 \textit{J. Pure Appl. Algebra}  \textbf{23}  (1982)  37-65.}

\bibitem{NN}{E. A. Navas and S. Nelson. On Symplectic Quandles. 
In Preparation.}

\bibitem{N}{S. Nelson. Classification of finite Alexander quandles.
\textit{Top. Proc.} \textbf{27} (2003) 245-258.}

\bibitem{NW}{S. Nelson and C-Y. Wong. On the orbit decompostion of finite 
quandles. \textit{J. Knot Theory Ramifications} \textbf{15} (2006) 761-772.}

\bibitem{Y}{D. Yetter. Quandles and monodromy. 
\textit{J. Knot Theory Ramifications} \textbf{12} (2003) 523-541.}

\end{thebibliography}
\end{document}